\documentclass[12pt]{article}

\usepackage{amsmath,amssymb,amsfonts,theorem,makeidx,latexsym,epsfig,,subfigure}

\newtheorem{defn}{Definition}[section]

\newtheorem{lemma}[defn]{Lemma}

{\theorembodyfont{\rmfamily}

\newtheorem{ex}[defn]{Example}}

\newtheorem{thm}[defn]{Theorem}

\newtheorem{prop}[defn]{Proposition}

\newtheorem{cor}[defn]{Corollary}

\newtheorem{rem}[defn]{Remark}

\numberwithin{equation}{section}

\newcommand{\h}{{\cal H}}

\newcommand{\afh}{ \forall f \in \h}

\newcommand{\mn}{\mathbb N}

\newcommand{\mc}{\mathbb C}

\def\bp{{\noindent\bf Proof. \ }}

\def\ep{\hfill$\square$\par\bigskip}

\def\bqs{\begin{equation}}

\def\eqs{\tag*{$\square$}\end{equation}\par\bigskip}

\def\la{\langle}

\def\ra{\rangle}

\def\ftk{\{f_k\}_{k=1}^\infty}

\def\ctk{\{c_k\}_{k=1}^\infty}

\def\gtk{\left\{g_k\right\}_{k=1}^\infty}

\def\htk{\{h_k\}_{k=1}^\infty}

\def\etk{\{e_k\}_{k=1}^\infty}

\def\suk{\sum_{k=1}^\infty}

\def\nl{\left|\left|}

\def\nr{\right|\right|}

\def\span{\overline{\text{span}}}

\def\Span{\text{span}}

\def\vn{\vspace{.1in}\noindent}

\def\bop{\begin{op}\rm}

\def\eop{\end{op}}

\def\bee{\begin{eqnarray}}

\def\ene{\end{eqnarray}}

\def\bes{\begin{eqnarray*}}

\def\ens{\end{eqnarray*}}

\def\bei{\begin{itemize}}

\def\eni{\end{itemize}}

\def\bt{\begin{thm}}

\def\et{\end{thm}}

\def\bc{\begin{cor}}

\def\ec{\end{cor}}

\def\bpr{\begin{prop}}

\def\epr{\end{prop}}

\def\bl{\begin{lemma}}

\def\el{\end{lemma}}

\def\bd{\begin{defn}}

\def\ed{\end{defn}}

\def\bex{\begin{ex}}

\def\enx{\end{ex}}

\def\bfi{\begin{fig}}

\def\efi{\end{fig}}

\def\k{{\cal K}}
\def\ptk{\{\psi_k\}_{k=1}^\infty}

\title{Completion versus removal of redundancy by perturbation}

\date{\today}

\author{Ole Christensen, Marzieh Hasannasab}

\begin{document}

\maketitle

\begin{abstract} A sequence $\gtk$ in a Hilbert space $\h$ has the expansion property if each $f\in \span \gtk$
has a representation $f=\suk c_k g_k$ for some scalar coefficients $c_k.$ In this paper we analyze the question
whether there exist small norm-perturbations of $\gtk$
which allow to represent all $f\in \h;$ the answer turns
out to be yes for frame sequences and Riesz sequences,
but no for general basic sequences. The insight gained from the analysis is used to address a somewhat dual
question, namely, whether it is possible to remove redundancy from a sequence with the expansion property via small norm-perturbations; we prove that the answer is yes for frames $\gtk$ such that $g_k\to 0$ as $k\to \infty,$ as well as for frames with finite excess. This
particular question is motivated by recent progress in
dynamical sampling.
\end{abstract}

\begin{minipage}{120mm}


\end{minipage}
\

\section{Introduction}

Let $\h$ denote a separable infinite-dimensional Hilbert space and suppose
that a given sequence $\gtk$ in $\h$ has the {\it expansion property,} i.e., that each $f\in \span \gtk$ has a representation
\bee \label{211702a} f= \suk c_k g_k\ene for certain coefficients
$c_k\in \mc.$ Our goal is to address the following question: when and how can we perform small norm-perturbations
on the sequence $\gtk$ and hereby obtain a sequence
$\ptk$ such that {\it arbitrary}
elements $f\in \h$ have an expansion
$  f= \suk c_k \psi_k$ for certain coefficients
$c_k\in \mc?$

Formulated as above, the question is clearly a
{\it completion  problem.}  We will show that the completion problem has an affirmative answer for the so-called Riesz
sequences and frame sequences, but not for general basic
sequences; along the way we
also consider a number of other completion problems.
Interestingly, the insight gained from the above analysis
can be used to address a somewhat dual question: when and how
can a {\it redundant} system $\gtk$ be turned into
a complete but nonredundant system $\ptk$ by small
norm-perturbations? We will provide a positive answer
to this question for a number of frames, in particular, for the so-called {\it near-Riesz bases}
introduced by Holub in \cite{Ho}. Additional
motivation for this particular question will be provided
at the end of the paper.

The paper is organized as follows. In the rest of the introduction we set the stage by providing a number of definitions and results from the literature. In Section \ref{211702c} we present the results about the completion problem; the dual problem concerning removal of redundancy is considered in Section \ref{211702d}.

A sequence $\gtk$ in
the Hilbert space $\h$ is called a {\it frame for $\h$} if there
exist constants $A,B>0$ such that
\bee \label{211902b} A \, ||f||^2 \le \suk | \la f, g_k\ra|^2\le B\, ||f||^2, \, \afh;\ene
suitable numbers $A,B$ are called {\it lower, resp. upper
frame bounds.} The sequence $\gtk$ is called a {\it Bessel sequence} if at least  the right-hand inequality in \eqref{211902b} holds.
A frame which is at the same time a
basis, is called a {\it Riesz basis.} Note that several other characterizations of frames and Riesz bases exist,
e.g., in terms of operator theory. For example, if
$\etk$ is a given orthonormal basis for $\h,$ frames
for $\h$ are precisely the sequences $\{Ve_k\}_{k=1}^\infty$ where $V: \h \to \h$ is a bounded
surjective operator;  Riesz bases correspond precisely
to the case where the operator $V$ also is injective.
Finally, a sequence $\gtk$ which is a frame
for the (sub)space $\k:= \span \gtk,$ is called
a {\it frame sequence;}  Riesz sequences are defined
in the analogue way.

One of the key reasons for the interest in frames is
that a frame has the expansion property: in fact,
given any frame $\gtk,$ there exists a so-called
{\it dual frame} $\ftk$ such that
\bes f= \suk \la f, f_k\ra g_k, \, \forall f\in \h.\ens
In general, the dual frame $\ftk$ is not unique: indeed,
the case where $\gtk$ is a Riesz basis is characterized
precisely by the existence of a unique dual.
We refer to
\cite{CB} for more information about frames and Riesz bases, also about their history
and applications.

The following Lemma collects a number of well-known results concerning
norm-perturbations of various sequences with the expansion property.

\bl \label{211702e} Consider two sequences $\gtk,
\htk$ in $\h,$ satisfying that
\bes \suk ||g_k-h_k||^2 < A,\ens for a value of $A$
as specified below. Then the
following holds:
\bei \item[(i)] {\bf \cite{C4}}  If
$\gtk$ is a frame for $\h$ with lower bound $A,$  then
$\htk$ is a frame for $\h.$
\item[(ii)] {\bf \cite{C9}} If
$\gtk$ is a Riesz sequence with lower bound $A,$ then
$\htk$ is a Riesz sequence; furthermore,
$\text{codim}( \span \gtk)= \text{codim}( \span \htk).$
\eni \el

 Alternative
norm-perturbation conditions are formulated in
\cite{CLZ}; however, they need that we have access to
information about a dual frame, which is not
the case in the current paper.
Note also that a number
of classical results about norm-perturbation (typically for
orthonormal sequences) are collected in \cite{Y}.
Observe that more general perturbation results are available in the literature, typically formulated in terms of certain
operators rather than norm-perturbations;
see, e.g., \cite{CB} and the references therein.


\section{Completion via norm-perturbation} \label{211702c}

Our main interest is to consider the completion problem for sequences $\gtk$ having the expansion property.
However, we first state a number of other completion
properties, some of which will be needed in later proofs.
Given any sequence $\gtk$ in $\h,$ we
define its {\it excess} ${\cal E}(\gtk)$ as the
maximal number of elements that can be removed
without changing the spanned space, i.e.,
\bee \label{201610a} {\cal E}(\gtk):= \max \sharp \{
J \subset \mn \, \big|  \, \span \{g_k\}_{k\in \mn \setminus J}= \span \gtk \}.\ene

\bpr \label{new} Let $\gtk$ be a sequence in $\h.$
Then the following hold:

\bei \item[(i)] If $\gtk$ is  not norm-bounded below,
there
exists a complete  sequence $\ptk$  such that
\bee \label{211902f} || g_k-\psi_k|| \to 0 \, \mbox{as } \, k\to \infty;\ene
\item[(ii)] If  ${\cal E}(\gtk) \ge \mbox{codim}(\span \gtk)$,
there
exists a complete  sequence $\ptk$  such that
\eqref{211902f} holds.
\item[(iii)] If $\gtk$ is convergent,
there
exists a complete  sequence $\ptk$  such that
\eqref{211902f} holds; in particular, $\ptk$
converges to the same limit as $\gtk.$
\eni In all the stated cases, given any $\delta >0,$
the sequence $\ptk$ can be
chosen such that
$|| g_k-\psi_k|| \le \delta$ for all $k\in \mn.$
\epr

\bp For the proof of (i),
given $\delta >0,$  choose a frame $\ftk$ for $\h$ such
that $||f_k|| \le \delta$ for all $k\in \mn$ and
$||f_k|| \to 0$ as $k\to \infty;$ for
example,  letting $\etk$ denote any
orthonormal basis, we can take
\bes \ftk = \{\delta\, e_1, \frac{\delta}{\sqrt{2}}\, e_2,  \frac{\delta}{\sqrt{2}}\, e_2, \frac{\delta}{\sqrt{3}}\, e_3,\frac{\delta}{\sqrt{3}}\, e_3,\frac{\delta}{\sqrt{3}}\, e_3, \dots\}.\ens
Denote the lower frame bound for the frame $\{ f_k \}_{k=1}^\infty$ by $A.$
Choose now a
subsequence $\{g_{k_n}\}_{n=1}^\infty $ of
$\gtk$ such that $ ||g_{k_n}||^2 \le \frac{3  A}{\pi^2 n^2}, \, n\in \mn;$ then
\bes \sum_{n=1}^\infty  || f_n- \left( f_n+g_{k_n}\right)||^2 =  \sum_{n=1}^\infty
||g_{k_n}||^2 \le \frac{A}{2}.\ens
Using Lemma \ref{211702e}(i), this
implies that $\{ f_n+g_{k_n}\}_{n=1}^\infty$
is a frame for $\h$ and hence complete. Thus, the sequence
$\ptk$ formed from $\gtk$ by replacing the
elements $\{g_{k_n}\}_{n=1}^\infty $ by
$\{f_n+g_{k_n}\}_{n=1}^\infty$
will satisfy the requirements.

For the proof of (ii) we first  assume
additionally that   $M:=\mbox{codim}(\span \gtk)$ is finite. Without loss of generality and only for notational convenience, assume that
the sequence $\gtk$ is ordered such that
$g_1, \dots, g_M\in \span\{g_k\}_{k=M+1}^\infty,$ and
take an orthonormal basis $\{e_k\}_{k=1}^M$
for the orthogonal complement $(\span \gtk)^\perp.$ Then the sequence
\bes \{\psi_k\}_{k=1}^\infty =
\{g_1+ \delta e_1, g_2+ \frac{\delta}{2}\, e_2,\dots, g_M+ \frac{\delta}{M} e_M,
g_{M+1}, g_{M+2}, \dots\}\ens
satisfies the requirements.
The case
where ${\cal E}(\gtk)= \mbox{codim}(\span \gtk)=\infty$
is similar and only requires minor notational
modifications.

For the proof of (iii), assume that
the sequence $\gtk$ converges to $f\in \h.$ Given $\delta >0,$
choose $K\in \mn$ such that $||f-g_k|| \le \delta/2$
for $k\ge K.$ Let $\etk$ denote an orthonormal basis
for $\h$ and define $\ptk$ by

\bes \psi_k:= \begin{cases}  g_k, &\mbox{if } \, k=1, \dots, K-1; \\
f, &\mbox{if } \, k=K; \\
f+ \frac{\delta}{2^{k-K}}e_{k-K}, &\mbox{if } \, k>K.\end{cases} \ens
Then $\Span \etk \subseteq \Span \ptk,$ so
$\Span \ptk$ is clearly complete. Furthermore, for
$k\ge K,$
\bes || g_k-\psi_k|| \le ||g_k-f|| + ||f-\psi_k|| \le\delta,\ens and $|| g_k-\psi_k||\to 0$ as
$k\to \infty.$
\ep

We are now ready to consider the completion problem
for Riesz sequences and frame sequences. The proofs
rely on an interesting result proved recently by
V. Olevskii:

\bl \label{Olevskii} {\bf \cite{VO1,VO2}}  If $\etk$ is an
orthonormal sequence in $\h,$ there
exists an orthonormal basis $\{ \chi_k\}_{k=1}^\infty$ for $\h$ such that
\bes || e_k-\chi_k || \to 0 \, \mbox{as } \, k\to \infty.\ens
In addition, given any $\delta >0,$
the sequence $\{ \chi_k\}_{k=1}^\infty$ can be
chosen  such that \\
$|| e_k-\chi_k || \le \delta$ for all $k\in \mn.$
\el

\newpage

\bt  \label{marzieh} Let $\gtk$ be a sequence in $\h.$
Then the following hold:

\bei
\item[(i)] If $\gtk$ is a Riesz sequence, there
exists a Riesz basis $\ptk$ for $\h$ such that
\bee \label{211902g} || g_k-\psi_k|| \to 0 \, \mbox{as } \, k\to \infty.\ene
\item[(ii)] If $\gtk$ is a frame
sequence,
there
exists a frame $\ptk$ for $\h$ such that
\eqref{211902g} holds.

\item[(iii)] If $\gtk$ is a Bessel
sequence,
there
exists a complete Bessel sequence $\ptk$  such that
\eqref{211902g} holds.
\eni
In all the stated cases, given any $\delta >0,$
the sequence $\ptk$ can be
chosen such that
$|| g_k-\psi_k|| \le \delta$ for all $k\in \mn.$
\et

\bp We first prove (iii). Thus, let $\gtk$ be a Bessel sequence in $\h,$
and let  $\k:= \span \gtk;$ we can assume that
$\k^\bot \neq \{0\}.$ Also, if $\k$ is finite-dimensional,
the results follow from Proposition \ref{new} (ii), so we will assume
that $\k$ is infinite-dimensional.
Now, by the standard properties
of a Bessel sequence \cite{CB}, choose an orthonormal basis
$\etk$ for $\k$ and a bounded operator
$U: \k \to \k$ such that $g_k=Ue_k, \, k\in \mn.$
Associated with the orthonormal
sequence $\etk$, choose the orthonormal basis
$\{ \chi_k\}_{k=1}^\infty$ for $\h$ as in Lemma \ref{Olevskii}, and define a bounded operator $V: \h \to \h$ by
\bee \label{200709a}  V=U \mbox{ on } \k, \, \,
V=I \mbox{ on } \k^\bot.\ene
Since the range of the operator $U$
contains the vectors $\gtk,$ it is dense
in $\k.$ Thus, the range of the operator $V$ is dense
in $\h;$ this implies that the sequence
$\ptk:=\{ V\chi_k\}_{k=1}^\infty$ is complete in $\h.$
A direct calculation reveals that $\ptk$
is a Bessel sequence. Furthermore, for all $k\in \mn,$
\bes ||g_k-\psi_k||= ||Ue_k - V \chi_k|| = ||Ve_k- V\chi_k||
\le ||V||\, || e_k- \chi_k||.\ens
Since the operator $V$ only depends on the sequence
$\gtk$ (and the fixed choice of $\etk$), this proves
(iii).  This also gives the proof
of (i)-(ii). Indeed, if $\gtk$ is a frame
sequence, the range of the operator $U$ equals
$\k,$ which implies that the range of the operator
$V$ equals $\h,$ and hence $\ptk$
is a frame for $\h;$ and if $\ftk$ is a Riesz sequence,
the operator $U:\k \to \k$ is bijective,
implying that $V:\h \to \h$ is bijective, and
hence that $\ptk$ is a Riesz
sequence. \ep

\begin{rem} \rm{Despite the fact that $\delta>0$
can be chosen arbitrarily small in Theorem
\ref{marzieh}, there is a restriction on
how ``close" the sequence $\ptk$ can be to the
sequence $\gtk$. Indeed, if $\gtk$ is a (non-complete) Riesz sequence with
lower bound $A,$ then the sequence $\ptk$ in Theorem
\ref{marzieh}(i) must satisfy that  \bee \label{212102a} \suk ||g_k-\psi_k||^2 \ge A;\ene otherwise
Lemma \ref{211702e}(ii) would imply that $\ptk$ is non-complete as well.
A similar result holds for frame sequences, although  the lower bound on the infinite sum in \eqref{212102a} will involve the gap between two particular subspaces
of $\h;$ see \cite{C9,CDL} for more detailed information.}\end{rem}

Theorem \ref{marzieh} makes it natural to ask whether
a basic sequence (i.e.,
a Schauder basis for a subspace) also can be extended
to a Schauder basis for $\h$ by small norm-perturbations
of the elements.  The following example shows that
the answer is no in general,
unless additional assumptions are added.

\bex \label{m2} Let $\etk$ denote an orthonormal basis
for $\h$ and consider the sequence
\bes \gtk= \{2e_2, 4e_4, 6e_6, \dots\}=\{2k e_{2k}\}_{k=1}^\infty.\ens
Clearly $\gtk$ is a basic sequence. Now, given any
$\delta \in ]0, 2\sqrt{6}\pi^{-1}[,$ consider a sequence
$\{\psi_k\}_{k=1}^\infty $ in $\h$ such that
$|| g_k- \psi_k|| \le \delta$ for all $k\in \mn.$ Then
\bes \sum_{k=1}^\infty || e_{2k} - \frac1{2k}\, \psi_k||^2
= \sum_{k=1}^\infty \frac1{4k^2}\, || g_k- \psi_k||^2
\le \frac{\pi^2 \delta^2}{24} <1.\ens
Since $\{e_{2k}\}_{k=1}^\infty$ forms a Riesz
sequence with lower bound $A=1,$ Lemma~\ref{211702e}(ii) implies
that $\{(2k)^{-1} \psi_k\}_{k=1}^\infty$ also forms
a Riesz sequence, spanning a space of the
same codimension as $\{e_{2k}\}_{k=1}^\infty;$
in particular, $\{ \psi_k\}_{k=1}^\infty$
can not be complete in $\h,$ and hence is not a
Schauder basis for $\h.$  \ep \enx

\section{Removal of redundancy via norm-perturbations} \label{211702d}

In this section the focus is on sequences $\gtk$
having the expansion property on the entire
underlying Hilbert space $\h.$ Such expansions might
be redundant, i.e., a given $f\in \h$ might have expansions $f= \suk c_k g_k$ for more than one choice
of the scalar coefficients $\ctk.$ A typical example of a redundant sequence is a frame $\gtk$ which is not a Riesz
basis. Our goal is to show that for certain frames $\gtk$ the redundancy can be removed via small norm-perturbations of
the vectors $g_k.$

Our first observation, stated next, does not even need
the frame assumption or any other expansion property.
\bt \label{211601a} Consider any sequence $\gtk$ in $\h$ such that
$g_k\to 0$ as $k\to \infty.$ Then, given any $\delta >0,$
there exists a Riesz basis $\ptk$ for $\h$ such that
\bes ||g_k- \psi_k || \le \delta, \, \forall k\in \mn.\ens\et

\bp  First, given any $\delta >0,$ choose $K\in \mn$
such that $||g_k|| < \delta/2$ for $k\ge K.$ We will
now construct $\ptk$ recursively, of the form
$\psi_k:= g_k + \varphi_k$ with the vectors $\varphi_k$
chosen as described next. First,
take $\varphi_1\in \h$ such that $|| \varphi_1|| \le \delta$ and
$\psi_1 \neq 0.$ Then choose $\varphi_2\in \h$ such that $|| \varphi_2|| \le \delta$ and
$\{\psi_1, \psi_2\}$ is linearly independent. Continuing
recursively, we finally choose $\varphi_K\in \h$ such that $|| \varphi_K|| \le \delta$ and
$\{\psi_1, \psi_2, \dots, \psi_K\}$ is linearly independent.
Then $\{\psi_1, \psi_2, \dots, \psi_K\}$ is a Riesz basis for
the subspace $V:=\Span\{\psi_1, \psi_2, \dots, \psi_K\}.$
Now choose an orthonormal basis $\etk$ for $V^\bot$
and define $\psi_k$ for $k>K$ by $\psi_k:= \frac{\delta}{2}e_k.$ Then $\ptk$ is a Riesz basis for
$\h$ and $||g_k-\psi_k|| \le \delta$ for all $k\in \mn.$ \ep

The result in Theorem \ref{211601a} immediately
applies to a number of well-known frames in the literature:

\bex \label{211902c} We state a number of examples of frames $\gtk$ such that
$g_k\to 0$ as $k\to \infty:$

\bei \item[(i)]
Given any orthonormal basis $\etk$ for $\h,$
the family
\bes \gtk:= \{e_1, \frac1{\sqrt{2}}\,e_2,\frac1{\sqrt{2}}\,e_2, \frac1{\sqrt{3}}\,e_3, \frac1{\sqrt{3}}\,e_3, \frac1{\sqrt{3}}\,e_3,  \dots \}\ens
is a frame for $\h.$ Clearly $g_k\to 0$ as $k\to \infty.$
Note that this particular frame was used in the proof
of Proposition \ref{new}.
\item[(ii)] Let again $\etk$ be an orthonormal basis for
$\h$ and fix any $\alpha \in ]0,1[.$ Let $\lambda_\ell:=
1- \alpha^{-\ell}$ for $\ell\in \mn,$ and define the vectors \bes g_k:= \sum_{\ell=1}^\infty \lambda_\ell^k
\sqrt{1- \lambda_\ell^2} e_{\ell}, \, k\in \mn.\ens
Then $\gtk$ is a frame (the so-called Carleson frame),
a result proved by Aldroubi et al. in \cite{A1,A2}. It is
easy to see that $g_k \to 0$ as $k\to \infty.$  Note
that $\gtk$ is heavily redundant:  it can be proved that
for any $N\in \mn,$ any subfamily $\{g_{Nk}\}_{k\in \mn}$
of $\gtk$ is a redundant frame as well. From this point
of view it is surprising that $\gtk$ can be
approximated by a Riesz basis, as stated in
Theorem \ref{211601a}.
\item[(iii)] More generally than (ii), it was proved in
\cite{CHP} that any redundant frame that can be
represented as an operator orbit $\gtk= \{T^k \varphi\}_{k=1}^\infty$ for a bounded operator $T:\h \to \h$ and some $\varphi \in \h$ will have the property that
$g_k\to 0$ as $k\to \infty.$ \ep
\eni
\enx

In order to reach the next result we need the following Lemma.
Recall that the {\it deficit} of a sequence $\gtk$ is
defined as the codimension of the vector space
$\span \gtk.$

\bl \label{212101a} Let $\etk$ be an orthonormal basis for $\h.$
Given any $\delta>0$ and any $N\in \mn,$ there exists
an orthonormal system $\{\varepsilon_k\}_{k=1}^\infty$
with deficit $N$ such that
$||e_k-\varepsilon_k||\le \delta$ for all $k\in \mn.$ \el

\bp Take any orthonormal system $\{\varphi_k\}_{k=1}^\infty$
with deficit $N,$ and choose via Lemma \ref{Olevskii}
an orthonormal basis $\{\chi_k\}_{k=1}^\infty$
for $\h$ such
that $||\varphi_k - \chi_k||\le \delta$ for all $k\in \mn.$
Then, choose the unitary operator $U: \h \to \h$
such that $e_k= U\chi_k,$ and let $\varepsilon_k:=U\varphi_k, k\in \mn.$  Then $\{\varepsilon_k\}_{k=1}^\infty$ is an orthonormal system
with deficit $N,$ and
$|| e_k-\varepsilon_k||=||U\chi_k-U\varphi_k||
= ||\chi_k-\varphi_k|| \le\delta$
for all $k\in \mn,$ as desired.\ep

\bt \label{211902a} Consider a frame of the form $\gtk= \{g_k\}_{k=1}^N \cup \{g_k\}_{k=N+1}^\infty,$ where
$N\in \mn$ and $\{g_k\}_{k=N+1}^\infty$ is a Riesz basis for $\h.$ Then, given any $\delta>0,$ there exists a Riesz basis $\ptk $ such that $||g_k-\psi_k||\le \delta$ for
all $k\in \mn.$ \et

\bp First, consider an orthonormal basis for $\h$
indexed as $\{e_k\}_{k=N+1}^\infty$ and choose
the bounded bijective operator $V:\h \to \h$ such that
$g_k=Ve_k$ for $k=N+1, N+2, \dots.$ Using Lemma \ref{212101a}, choose an orthonormal system $\{\varepsilon_k\}_{k=N+1}^\infty$
with deficit $N$ such that
$||e_k-\varepsilon_k||\le \delta/{||V||}$ for  $k=N+1, N+2, \dots.$
Then, letting $\psi_k:= V\varepsilon_k, k=N+1, N+2, \dots,$ the
family $\{\psi_k\}_{k=N+1}^\infty$ is a Riesz sequence with deficit $N,$  and $||g_k-\psi_k|| = || Ve_k-V\varepsilon_k||
\le  \delta$ for $k=N+1, N+2, \dots.$

Now, consider the vector $g_N.$  If $g_N \notin
\span  \{\psi_k\}_{k=N+1}^\infty,$ let $\psi_N:=g_N;$
then $\{\psi_k\}_{k=N}^\infty$ is a Riesz sequence with deficit $N-1.$  On the other hand, if
$g_N \in
\span  \{\psi_k\}_{k=N+1}^\infty,$ choose any normalized
vector $\varphi_N \notin \span  \{\psi_k\}_{k=N+1}^\infty,$
and let $\psi_N:= g_N+  \delta \varphi_N;$
then again $\{\psi_k\}_{k=N}^\infty$ is a Riesz sequence with deficit $N-1,$ and $||g_k-\psi_k||
\le  \delta$ for $k=N,N+1, N+2, \dots.$
Applying now the same procedure on
$g_{N-1},g_{N-2}, \dots,g_1,$ we arrive at the desired
Riesz basis $\ptk$ in a finite number of steps. \ep

Interestingly, frames of the type considered in Theorem \ref{211902a} were called {\it near-Riesz bases } by
Holub in the paper \cite{Ho}; the above result provides
an additional reason for this name being very appropriate.

\begin{rem}  \rm{Despite the fact that $\delta>0$ can
be chosen arbitrarily small in Theorem \ref{211902a},
the Riesz basis $\ptk$
 must satisfy that
$\suk ||g_k-\psi_k||^2 \ge A,$ where $A$ is the lower
frame bound for $\gtk;$ otherwise
the results in \cite{CC1} show that $\ptk$ would
be a frame with the same excess as $\gtk.$}
\end{rem}

We want to point out that the proof of
Theorem \ref{211902a}  somewhat hides the
fact that it is highly nontrivial to get direct access to
the Riesz basis $\ptk,$ especially due to the intriguing and
deep construction by V. Olevskii  playing a key role
in the argument. The next example illustrates this
by a concrete construction.

\bex \label{212102b} Let  again $\etk$ be an orthonormal basis for $\h$ and consider the frame
\bes \gtk:= \{e_1, e_1, e_2, e_3, e_4, \dots\},
\ens
consisting of the orthonormal basis and
a single extra element. A natural way to try to remove
the redundancy would be to fix a small $\epsilon>0$
and let $\psi_1:=e_1$
and for $k>1, \psi_k:= \frac12 e_{k-1} + ( \frac12 + \epsilon) e_k.$ Then for any finite sequence
$\{c_k\}_{k=2}^\infty,$
\bes  \nl \sum_{k=2}^\infty c_k( (\frac12 + \epsilon)e_k-\psi_k) \nr^2= \frac14 \nl \sum_{k=2}^\infty c_ke_{k-1} \nr^2 = \frac14 \sum_{k=2}^\infty |c_k|^2.\ens
Observe that $\{e_1\}\cup \{(\frac12 + \epsilon)e_k\}_{k=2}^\infty$ is a Riesz basis
with lower bound  $\frac12 + \epsilon.$  Considering
$\ptk$ as a perturbation of this Riesz basis, it
now follows from the results in \cite{CC1}
that $\ptk$ is a Riesz basis for $\h.$ Note that
\bes || g_k-\psi_k||= \sqrt{\frac 14 + (\frac12  + \epsilon)^2};\ens however, this construction does
not allow us to obtain $|| g_k - \psi_k|| \le \delta$
when  $\delta < 2^{-1/2}\approx 0.7.$  In fact, in
order to obtain the result in Theorem \ref{211902a} for smaller values of $\delta,$ it would
be necessary to consider much more
complicated perturbations   $\ptk$
of $\gtk,$ making it highly nontrivial to do this in practice.
\ep \enx

\begin{rem}
\rm{The question of removal of redundancy is partly motivated by the research topic
{\it dynamical sampling,} introduced in the papers
\cite{AD,A1}. One of the key issues in dynamical sampling
is the construction of frames as orbits $\{T^k \varphi\}_{k=0}^\infty$ of a bounded operator $T:\h \to \h,$ for some $\varphi \in \h;$  we encountered such frames already in Example \ref{211902c}(ii) \& (iii).
Unfortunately it is very difficult to construct  such frames, and the only concrete examples available in the
literature are indeed Riesz bases
\cite{CH10} and the Carleson frame \cite{A1} considered
in Example \ref{211902c} (ii). Also, it was proved
in \cite{CH10} that a near-Riesz basis never has this
property. This raises the natural question whether a
near-Riesz basis can be approximated by a Riesz basis,
and hence by an orbit of a bounded operator;
Theorem \ref{211902a} confirms that this indeed is possible. We will phrase this consequence of
Theorem \ref{211902a} as a separate result, where we
index the given near-Riesz basis by $\{g_k\}_{k=0}^\infty$ for notational convenience:
}   \end{rem}

\bc \label{210103a} Consider any near-Riesz basis $\{g_k\}_{k=0}^\infty.$
Then, given any $\delta>0,$ there exists
$\varphi \in \h$ and a bounded operator $T:\h \to \h$
such that
\bes ||g_k-T^k \varphi|| \le \delta, \, \forall k\in \mn_0.\ens \ec

The results in Theorem \ref{211601a} and  Theorem \ref{211902a} do not cover
the standard (regular)
redundant Gabor frames and wavelet frames: they
consist of vectors with equal norm, and they have
infinite excess \cite{B}. Due to the complications
discussed in Example \ref{211601a} and the preceding
text it seems to be very difficult
to answer the question whether all frames indeed can be
approximated by a Riesz basis. At least for Gabor frames
and wavelet frames  we
can apply the following adaption
of the Feichtinger Theorem (finally proved in one
of its equivalent formulations in \cite{MSS}), showing
that any frame which is norm-bounded below can be approximated by a {\it finite}
collection of Riesz bases:

\bt \label{210103b} Let $\gtk$ be a frame which is norm-bounded below.
Then there exists a finite partition $\gtk = \bigcup_{j=1}^J
\{g_k\}_{k\in I_j}$ with the property that for each $\delta >0$
there exist Riesz bases $\{\psi_k\}_{k\in I_j}, j=1, \dots, J,$ for
$\h$ such that $||g_k-\psi_k||\le \delta$ for all $k\in \mn.$ \et

\bp Choose according to the Feichtinger Theorem
a finite partition  $\gtk = \bigcup_{j=1}^J
\{g_k\}_{k\in I_j}$ such that each sequence
$\{g_k\}_{k\in I_j}, j=1, \dots,J$ is a Riesz sequence; using
Theorem 2.1 in \cite{CH10} we can shuffle the
elements around to ensure that each of the
index sets $I_j$ is infinite. Now the result follows
directly from Theorem \ref{marzieh}(i). \ep

The result in Theorem \ref{210103b} can be formulated
as a operator-theoretic result, similarly to
Corollary \ref{210103a}; we leave the precise formulation to the interested reader.

{\bf \vspace{.1in}

\noindent Ole Christensen\\
DTU Compute\\
Technical University of Denmark\\
Building 303 \\
2800 Lyngby  \\
Denmark \\
 Email: ochr@dtu.dk

\vn Marzieh Hasannasab \\
Institut f\"ur Mathematik\\
TU Berlin\\
Stra\ss e des 17. Juni 136
\\ 10623 Berlin \\ Germany
\\
Email: hasannas@math.tu-berlin.de

}

\end{document}